\documentclass[11pt]{article}
\usepackage{amsmath, amsfonts, amsthm, amssymb,verbatim,graphicx,color}

\numberwithin{equation}{section}

\newcommand \la \langle
\newcommand \ra \rangle

\newcommand{\ekp}{e_{K}^{+}}
\newcommand{\ekm}{e_{K}^{-}}

\newcommand \TT     {\mathcal{T}}
\newcommand \Tcal   {\TT}

\newcommand \del        \partial

\newcommand \be     {\begin{equation}}
\newcommand \ee     {\end{equation}}
\newcommand \RR      {\mathbb{R}}
\newcommand \eps     \epsilon

\let\oldmarginpar\marginpar
\renewcommand\marginpar[1]{\-\oldmarginpar[\raggedleft\footnotesize #1]%
{\raggedright\footnotesize #1}}
\begin{document}

\title{Finite Volume Method for the Relativistic Burgers Model on a (1+1)-Dimensional de Sitter Spacetime}

\author{
Tuba Ceylan$^1$   and Baver Okutmustur$^2$}

\date{}

\maketitle

\footnotetext[1] {Central Bank of the Republic of Turkey (TCMB). E-mail: {\sl ceylanntuba@gmail.com}}
\footnotetext[2]{Department of Mathematics, Middle East Technical University (METU), 06800 Ankara, Turkey.
E-mail : {\sl  baver@metu.edu.tr}
\
\newline
\textit{{\bf Mathematics Subject Classification.}} 35L65, 35Q53, 65M08, 76N10, 83A05
\newline
\textit{{\bf Key Words and Phrases.}} relativistic Burgers equation; Euler system; de  Sitter metric; de  Sitter backgrounds; \mbox{finite volume method}; Godunov scheme}


\begin{abstract}
Several generalizations of the relativistic models of Burgers equations have recently been  established  and developed  on  different spacetime geometries. In this work, we take into account  the de Sitter  spacetime geometry, introduce   our relativistic  model by a technique based on the vanishing pressure Euler equations of relativistic compressible fluids on a (1+1)-dimensional background and construct a second order Godunov type finite volume scheme to examine numerical experiments within an analysis of the cosmological constant.  Numerical  results demonstrate the efficiency of the method for solutions containing shock  and rarefaction waves.

\end{abstract}

\section{Introduction}

A fundamental relativistic Burgers equation has recently been derived and generalized to different spacetime geometries by LeFloch and his collaborators \cite{CLO, CO, CO2,  LMO}. The first model has been obtained by establishing a hyperbolic balance law  satisfying the  Lorentz invariance property  for flat geometry  and then its relativistic generalizations were extended to the Schwarzschild,  Friedmann--Lemaitre--Robertson--Walker  (FLRW), de Sitter and Schwarzschild--de Sitter   \mbox{spacetimes \cite{Ceylan, CLO, CO, CO2, LMO}.}

In this article, we consider  the ``de Sitter'' (dS) spacetime which is a member  of the family of the FLRW geometry.  The line elements of the FLRW and the dS geometries share certain  common properties. In particular, these metrics  are solutions to the Einstein field equations. Furthermore, the Minkowski metric can be derived from particular cases of both the FLRW and the dS metrics.  The inspiration of derivation of the  relativistic Burgers equation on the dS spacetimes is based on the common  properties of these geometries.  We follow the instruction via \cite{CLO, LMO, LM} and take into account  the relativistic Euler equations on a given curved background $M$ by the following general formulation of balance laws:

\be
\label{Euler0}
\aligned
&\nabla_\alpha \big( T^{\alpha\beta} (\rho,u)\big) =0,\\
&T^{\alpha\beta} (\rho,u)= \rho c^2 u^{\alpha}u^{\beta} + p(\rho) \big( u^{\alpha}u^{\beta}  +g^{\alpha \beta}\big),
\endaligned
\ee
where $ \nabla$ is the covariant derivative operator, $T^{\alpha\beta}$ is the energy-momentum tensor of perfect fluids,  $\rho \geq 0$ denotes the mass-energy density of the fluid, $c$ is the light speed and the unit timelike vector field $u^{\alpha}$ represents the  velocity of the fluid so  that $g_{\alpha \beta} u^{\alpha} u^{\beta}=-1$.

\subsection{Motivation of the Paper}

One of the basic nonlinear hyperbolic conservation laws:
\be
\label{B}
\frac{\del v}{\del t} + \frac {\del}{\del x}\big(\frac {\,\,v^2}{2}\big)=0, \, v=v(t,x), \, t>0, x\in \RR
\ee
 is the (inviscid) Burgers equation which has various applications in many physical fields. Indeed, this equation  is a simplified version of the Euler equations of compressible fluids
\be
\label{EE1}
\aligned
&\frac{\del \rho}{\del t} + \frac{\del}{\del x} (\rho  v)= 0,\\
&\frac{\del }{\del t}(\rho v) + \frac{\del}{\del x}( \rho v^2 + p(\rho) ) = 0,
\endaligned
\ee
where $v$ denotes the velocity, $\rho  $  the density,  and $p$ the pressure of the fluid. We recover the inviscid Burgers equation from this system by first imposing that  $p \equiv 0$ in Equation \eqref{EE1}. Then, we take partial derivatives for both equations, combine  them  and finally   deduce the Equation \eqref{B}. This analysis was generalized for  the relativistic models in  the articles \cite{CLO, CO, LMO} in order  to get the relativistic Burgers equations on the related geometries.

In the present article, we apply this analysis by considering the relativistic Euler equations on a smooth, time-oriented, curved (1+1)-dimensional { dS}  {  spacetime}  in order to derive the relativistic Burgers model.
We pursue and cite  the article by Ceylan and Okutmustur \cite{CO} for the details of derivation of  the model on a dS  background. Once we obtain our model, we then construct a second order Godunov type finite volume scheme  in order to inspect the stability and efficiency of the method.
 Following the papers by LeFloch \textit{et al.} \cite{CLO, LMO},  the relativistic Burgers equations are derived from  Equation \eqref{Euler0} on a (1+1)-dimensional background  by imposing a vanishing pressure to Equation \eqref{Euler0}. The derived equations satisfy the Lorentz invariance property and in the case where $c\to \infty$,  the classical (non--relativistic) Burgers equation is recovered directly from the given system. In paper \cite{ LMO}, the curved spacetime is considered to be flat and then generalized to a Schwarzschild background; whereas in  \cite{ CLO}, an FLRW background is taken into account. In both  \cite {CLO, LMO}, the derived relativistic models are used for the numerical experiments via a second-order accurate finite volume method  based on the geometric formulation Equation \eqref{Euler0}. A similar work is attained by Ceylan and Okutmustur for the Schwarzschild--de Sitter geometry \cite{CO2}. For  the convergence  and geometric formulation of finite volume schemes  on Lorentzian (curved) manifolds, we follow the articles by  LeFloch \textit{et al.} \cite {ALO,  LO}. The  theory of conservation laws on the flat space is introduced by Kruzkov \cite{SNK}. We refer the reader to the following work by Groah, Smoller and Temple \cite{GST} for the shock wave interactions and further details in general relativity. On the other hand, we cite a recent work on balance laws posed on (1+1)-dimensional spacetimes by Gosse \cite{Gosse} where a  locally inertial Godunov scheme is presented numerically.
  Moreover,  further details  of second order Godunov schemes are presented in the book of Guinot \cite{VG} and the article by Van Leer \cite{L}. The reader can find a more general instruction for the theory of general relativity and related concepts about Einstein's theory in the book of Wald \cite{W}.

An outline of this paper is as follows. We  introduce the basic features of a  dS geometry and its  metric in the first part. We then derive the  Euler equations by means of the Christoffel symbol where we refer the reader to the article \cite{CO} for  further details. The next issue is deriving the particular cases of the relativistic Burgers equations  depending on the cosmological constant parameter $\Lambda$. \linebreak Then, we present a geometric formulation of  finite volume approximation on a curved background which follows by a construction of finite volume scheme on local coordinates.  We then establish a Godunov type second order scheme for the given model. We also compare the non-relativistic  and  relativistic Burgers equations depending on  different values of  $\Lambda$.  The effects of the $\Lambda$ on the numerical scheme are studied in the numerical experiments part. We finish the article by investigating the efficiency of the scheme with different $\Lambda$ values. The results demonstrate that the scheme is consistent with the conservative form of our model which yields correct computations of weak solutions containing shock and rarefaction waves.

\section{The dS Metric and the Euler System }

The dS spacetime is a particular case of the Lorentzian manifold and its metric provides  a cosmological solution to the  Einstein field equations similar to the FLRW metric. The line element (metric) of the dS background involves a constant so called ``cosmological constant'' and it is denoted by $\Lambda$.  Einstein firstly introduced this cosmological constant in his field equations. Depending on the  sign of the cosmological constant $\Lambda$, we distinguish the de Sitter and the Anti-de Sitter spacetimes; namely, if $\Lambda >0$, the geometry is called the ``de Sitter spacetime'' and, if $\Lambda <0$, it is called the ``Anti-de Sitter spacetime''. Particularly, if $\Lambda=0$, the metric turns to be a Minkowski metric and thus we have a flat geometry.

The corresponding line element for a (3+1)-dimensional case in terms of  time $t$, the radial $r$ and angular coordinates  $\theta$ and $\varphi$  is given by:
\be
\label{flrw}
g =-(1-\Lambda r^2)dt^2+{1 \over 1-\Lambda r^2}dr^2+r^2(d\theta^2+\sin^2\theta d\varphi^2).
\ee
Thus, Equation \eqref{flrw} is called the ``de Sitter metric'' if $\Lambda >0$ and the ``Anti-de Sitter metric'' if $\Lambda<0$.
Notice that  Equation \eqref{flrw} can be rewritten in matrix form as
$$
(g_{ij})= \begin{pmatrix}
-(1-\Lambda r^2) & 0 & 0 & 0  \\
 0 & {1 \over 1-\Lambda r^2} & 0 & 0 \\
 0 & 0& r^2 & 0 \\
 0 & 0& 0 & r^2\sin^2\theta
\end{pmatrix}
$$
with  the non-zero ``covariant'' elements
\begin{eqnarray*}
\label{g1}
&g_{00}=-(1-\Lambda r^2), \quad  g_{11}={1 \over 1-\Lambda r^2}\\
 &g_{22}=r^2, \qquad \qquad \,  g_{33}= r^2\sin^2\theta
\end{eqnarray*}
and the inverse of $(g_{ij})$ is
$$
(g^{ij})= \begin{pmatrix}
{1 \over \Lambda r^2-1} & 0 & 0 & 0  \\
 0 & {1-\Lambda r^2} & 0 & 0 \\
 0 & 0& {1\over r^2} & 0 \\
 0 & 0& 0 & {1\over r^2\sin^2\theta}
\end{pmatrix}
$$
with  the corresponding  ``contravariant'' elements
\begin{eqnarray*}
\label{g2}
&g^{00}={1 \over \Lambda r^2-1}, \quad \quad g^{11}=1-\Lambda r^2  \\
&g^{22}={1\over r^2}, \qquad \qquad g^{33}= {1\over r^2\sin^2\theta}.
\end{eqnarray*}
Here,
$$
g^{ik}g_{kj}=\delta_j^i=\begin{cases} 1 &\mbox{if } i = j\\
0 & \mbox{if } i\neq j \end{cases}
$$
is satisfied where $\delta_j^i$ is the Kronecker's delta function.

\subsection{Christoffel Symbols for the dS Spacetime}

The Christoffel symbols are defined by
\be
\label{CRF}
\Gamma_{\alpha\beta}^\mu={1\over 2}g^{\mu \nu}(-\del_\nu  g_{\alpha\beta}+ \del_\beta g_{\alpha\nu}+ \del_\alpha g_{\beta\nu}),
\ee
where the terms $\alpha, \beta, \mu, \nu \in \{0,1,2, 3\}$.
By substituting  $\alpha, \beta, \mu, \nu \in \{0,1,2, 3\}$ in  Equation \eqref{CRF}, we calculate all the terms  of $\Gamma_{\alpha\beta}^\mu$ in the following:
$$
\aligned
&\Gamma_{00}^0=0,\\
&\Gamma_{01}^0={\Lambda r \over \Lambda r^2-1},\\
&\Gamma_{01}^0=\Gamma_{10}^0={\Lambda r \over \Lambda r^2-1},\\
&\Gamma_{11}^1={\Lambda r \over 1- \Lambda r^2},\\
&\Gamma_{00}^1=\Lambda r( \Lambda r^2-1),\\
&\Gamma_{22}^1= r( \Lambda r^2-1),\\
&\Gamma_{33}^1= r( \Lambda r^2-1)\sin^2\theta,\\
&\Gamma_{12}^2=\Gamma_{21}^2=\Gamma_{13}^3=\Gamma_{31}^3={1 \over r},\\
&\Gamma_{33}^2=- \sin\theta \cos\theta,\\
&\Gamma_{23}^3=\Gamma_{32}^3=\cot \theta.\\
\endaligned
$$

\subsection{Energy-Momentum Tensor}
%

We take into account  our spacetime to be (1+1)-dimensional so that
the  solutions to the Euler system depend only on the time
variable $t$ and radial variable $r$ and   the angular components $\theta, \varphi$ vanish. It follows that
\be
\label{unit}
(u^\alpha) = (u^0, u^1, 0,0), \mbox { with } u^{\alpha} u_{\alpha}=-1
\ee
and
$$
\aligned
u^{\alpha} u_{\alpha}&=u^0 u_0+u^1 u_1,\\
&=g_{00}(u^0) (u^0)+g_{11}(u^1) (u^1),\\
&= g_{00}(u^0)^2+g_{11}(u^1)^2,\\
&=-1,
\endaligned
$$
where  $(u^0)^2$ and $(u^1)^2$ are unitary vectors. Next, we substitute the covariant terms $g_{00}$ and $g_{11}$
in the above relation which yields
\be
\label{unit}
-1 = -(1-\Lambda r^2) (u^0)^2+{1\over 1-\Lambda r^2}(u^1)^2.
\ee

Next, we define a velocity component $v$ depending on $u^0$ and $u^1$ as follows:
\be
\label{velocity}
v :=  {c \over (1-\Lambda r^2)} {u^1 \over u^0}.
\ee
Examining the relations \eqref{unit} and \eqref{velocity} together, it follows that
\be
\label{u}
\aligned
(u^0)^2 = {c^2 \over  (1-\Lambda r^2)(c^2-v^2) },\quad
\textnormal {and}~(u^1)^2 = {v^2  (1-\Lambda r^2) \over(c^2-v^2) }.
\endaligned
\ee
It remains to substitute the Equation \eqref{u} terms into the energy momentum tensor for perfect fluids equation, which is defined by
\be
\label{tensors}
T^{\alpha\beta} = (\rho c^2 + p) \, u^\alpha \, u^\beta + p \, g^{\alpha\beta}.
\ee
Thus, by using Equations \eqref{u} and \eqref{tensors}, the tensors are evaluated as follows:
$$
\aligned
&T^{00} ={ \rho c^4+pv^2 \over (c^2-v^2)(1-\Lambda r^2)},\\
&T^{01} = T^{10} = {cv( \rho c^2+p) \over (c^2-v^2)},\\
&T^{11} = {c^2(1-\Lambda r^2)(v^2 \rho +p) \over (c^2-v^2)},\\
&T^{22} = {p \over r^2},\\
&T^{33} = {p \over r^2 \sin^2\theta},
\endaligned
$$
$$
T^{02} =T^{03}=T^{12}=T^{13} =T^{20} =T^{21}
 =T^{23} =  T^{30}=T^{31}=T^{32}=0.
$$

\subsection{The Euler System on a $(1+1)$-Dimensional dS Background with Zero Pressure}

We are about  to write the Euler equations and deduce the presureless Euler system from these equations.  We recall that the Euler system \eqref{Euler0} is given by
$$\nabla_{\alpha} T^{\alpha \beta}=0,
$$
which is equivalent to
\be
\label{Euler000}
\del_\alpha T^{\alpha\beta} +\Gamma_{\alpha\gamma}^\alpha  T^{\gamma\beta}+\Gamma_{\alpha\gamma}^\beta  T^{\alpha\gamma}=0,
\ee
where $T^{\alpha \beta}$ terms are calculated in the previous subsection.  In order to derive the desired system, we  start by substituting  $\beta=0$ in Equation \eqref{Euler000} to get
\be
\label{Tens1}
\del_\alpha T^{\alpha 0} +\Gamma_{\alpha\gamma}^\alpha  T^{\gamma 0}+\Gamma_{\alpha\gamma}^0  T^{\alpha\gamma}=0.
\ee
Then, we put $\alpha, \gamma \in \{0,1,2,3\}$ into the above equation which yields
$$
\aligned
&\del_0 T^{0 0}+\Gamma_{0 0}^0  T^{0 0}+\Gamma_{0 0}^0  T^{0 0}+\Gamma_{0 1}^0  T^{1 0}+\Gamma_{0 1}^0  T^{0 1}+\Gamma_{0 2}^0  T^{2 0}+\Gamma_{0 2}^0  T^{0 2} +\Gamma_{0 3}^0  T^{3 0}\\&
+\Gamma_{0 3}^0  T^{0 3}
+\del_1 T^{1 0} + \Gamma_{1 0}^1  T^{0 0}+ \Gamma_{1 0}^0  T^{1 0}+ \Gamma_{1 1}^1  T^{1 0}+ \Gamma_{1 1}^0  T^{11}+ \Gamma_{1 2}^1  T^{2 0}+ \Gamma_{1 2}^0  T^{12}\\& + \Gamma_{1 3}^1  T^{3 0}+ \Gamma_{1 3}^0  T^{13} +\del_2 T^{2 0}+ \Gamma_{20}^2  T^{0 0}+ \Gamma_{2 0}^0  T^{2 0}+ \Gamma_{2 1}^2  T^{1 0}+ \Gamma_{2 1}^0  T^{2 1}+ \Gamma_{22}^2  T^{2 0}\\& + \Gamma_{22}^0  T^{22}+ \Gamma_{23}^2  T^{3 0}+ \Gamma_{23}^0  T^{23}
+\del_3 T^{3 0}+ \Gamma_{3 0}^3  T^{0 0}+ \Gamma_{3 0}^0  T^{3 0}+ \Gamma_{3 1}^3  T^{1 0}+ \Gamma_{3 1}^0  T^{31}\\& + \Gamma_{32 }^3  T^{2 0}+ \Gamma_{32}^0  T^{32}+\Gamma_{33}^3  T^{3 0}+ \Gamma_{33}^0  T^{33}=0.
\endaligned
$$
Next, for $\beta=1$ in Equation \eqref{Euler000}, we have
\be
\label{Tens2}
\del_\alpha T^{\alpha 1} +\Gamma_{\alpha\gamma}^\alpha  T^{\gamma 1}+\Gamma_{\alpha\gamma}^1  T^{\alpha\gamma}=0,
\ee
and putting  $\alpha, \gamma \in \{0,1,2,3\}$ into the above equation yields
$$
\aligned
&\del_0 T^{0 1}+\Gamma_{0 0}^0  T^{0 1}+\Gamma_{0 0}^1  T^{0 0}+\Gamma_{0 1}^0  T^{1 1}+\Gamma_{0 1}^1  T^{0 1}+\Gamma_{0 2}^0  T^{2 1}+\Gamma_{0 2}^1  T^{0 2}
+\Gamma_{0 3}^0  T^{3 1}\\& +\Gamma_{0 3}^1  T^{0 3}
+\del_1 T^{1 1} + \Gamma_{1 0}^1  T^{0 1}+ \Gamma_{1 0}^1  T^{1 0}+ \Gamma_{1 1}^1  T^{1 1}+ \Gamma_{1 1}^1  T^{11}
+ \Gamma_{1 2}^1  T^{2 1}+ \Gamma_{1 2}^1  T^{12}\\& + \Gamma_{1 3}^1  T^{3 1}+ \Gamma_{1 3}^1  T^{13} +\del_2 T^{2 1}+ \Gamma_{20}^2  T^{0 1}+ \Gamma_{2 0}^1  T^{2 0}
+ \Gamma_{2 1}^2  T^{1 1}+ \Gamma_{2 1}^1  T^{2 1}+ \Gamma_{22}^2  T^{2 1}\\&+ \Gamma_{22}^1  T^{22}+ \Gamma_{23}^2  T^{3 1}+ \Gamma_{23}^1  T^{23}+\del_3 T^{3 1}
+ \Gamma_{3 0}^3  T^{0 1}+ \Gamma_{3 0}^1  T^{3 0}+ \Gamma_{3 1}^3  T^{1 1}+ \Gamma_{3 1}^1  T^{31}\\&+ \Gamma_{32 }^3  T^{2 1}+ \Gamma_{32}^1  T^{32}+\Gamma_{33}^3  T^{3 1}+ \Gamma_{33}^1  T^{33}
=0.
\endaligned
$$
Replacing the Christoffel symbols and the tensor terms with their calculated values and imposing zero pressure to this system yields
\be
\aligned
\label{EE2}
&\del_0\Big({c \over (1-\Lambda r^2) (c^2-v^2)}\Big)+\del_1 \big({v\over  {c^2-v^2}}\big)+{2v\over  r(c^2-v^2)}+{2 \Lambda r v\over  (\Lambda r^2 - 1)(c^2-v^2)}=0,\\
&\del_0\Big({cv \over c^2-v^2}\Big)+\del_1 \Big({v^2  (1-\Lambda r^2)\over  {c^2-v^2}}\Big)-{\Lambda r c^2 \over{c^2-v^2}}+{\Lambda r v^2 \over{c^2-v^2}}+{2 (1-\Lambda r^2)v^2 \over {r(c^2-v^2)}}=0,
\endaligned
\ee
which is  the ``Euler system with zero pressure'' on a (1+1)-dimensional dS background.
\section{Relativistic Burgers Equation on a (1+1)- Dimensional dS  Spacetime}

\subsection{Derivation of the Model}

We observed in the introduction part that the classical Burgers equation is deduced from the Euler equations  by imposing zero  pressure to the given system. In this section, we follow the same technique   in order to derive the relativistic Burgers equation on the dS  geometry. The investigation of  static and spatially homogeneous solutions of the given model is also an issue of this section.

We  start by rewriting  the first and second equations of  \eqref{EE2} with the  notation $\del_0=\del_t, \, \del_1=\del_r$ and combining these two equations (by taking partial derivatives of each term of \eqref{EE2} and summing these equations after elimination  of the vanishing terms) in order to derive the following single equation:
\be
\label{maineq}
\del_t v + (1-\Lambda r^2) \del_r({v^2\over 2})+\Lambda r(v^2-c^2)=0,
\ee
which is the desired so called ``relativistic Burgers equation on the dS background''. Notice that the left-hand side of the model can also be written in the conservative form as follows:
$$
\del_t v+\del_r\left((1-\Lambda r^2 ){v^2\over 2}\right)=\Lambda r(c^2-2v^2).
$$
\subsection{Limiting Case of the Derived Model}

The Equation \eqref{maineq} may differ depending on different values of $\Lambda$.  We recall that  the particular case of the dS metrics for the value $\Lambda=0$ gives the Minkowski metric. If we substitute $\Lambda=0$ in   the relativistic  Burgers equation \eqref{maineq}, we recover
the classical (inviscid) Burgers equation:
$$
\del_t v + \del_r({v^2\over 2})=0.
$$
In other words,  the limiting case of  the relativistic Burgers Equation \eqref{maineq} on a dS  spacetime yields the classical (non-relativistic) Burgers equation. This common property is also shared by the other relativistic models on the curved spacetimes (such as flat, Schwarzschild and FLRW).

\subsection{Static and Spatially Homogeneous Solutions for the Model}

In this subsection, we search for the static and homogeneous solutions of our model if they exist. To examine this, we need to write the conservative form of the derived model  (\ref{maineq}),  which is
\be
\label{sitter2}
\del_t v+\del_r\left((1-\Lambda r^2 ){v^2\over 2}\right)=\Lambda r(c^2-2v^2).
\ee
We start by seeking the static solutions, that is the $t$-independent solutions, to Equation \eqref{sitter2}.
Due to the $t$-independency, the term  $\del_t v$  vanishes,  and we obtain
\be
\label{static1}\del_r\left((1-\Lambda r^2 ){v^2\over 2}\right)=\Lambda r(c^2-2v^2).
\ee
Applying  the  change of variable
\be
\label{change}
\aligned
K=1-\Lambda r^2, \quad
L=c^2-v^2
\endaligned
\ee
 results in
$$
L=K\,N,
$$
where $N \in (0, c)$ is a constant parameter. It follows that
$$
c^2-v^2=N \,(1-\Lambda r^2).
$$
Hence, we find the ``static solutions'' described by
\be
\label{static2}
v_{static}= \pm \sqrt{c^2-N(1-\Lambda r^2)}.
\ee

On the other hand, we follow an analogue process for  the homogeneous solutions, that is the $r$-dependent solutions to the Equation \eqref{sitter2}. It is obvious to conclude that  there is no spatially homogeneous  ($r$-independent) solution since the equation
$$
\del_t v=\Lambda r(c^2-2v^2)
$$
depends already on $r$.

As a result, we conclude that  our model  (\ref{maineq})  admits only ``static solutions''.
\section{Finite Volume Method on a Curved Spacetime }

The finite volume method  is an important discretization technique for partial differential equations. In the present work, we use a  finite volume approximation for general balance laws of hyperbolic partial differential equations on an ($n$+1)-dimensional manifold following the\linebreak papers \cite{ALO, DAF, LMO, RL, LO}. We  consider the geometry in ($n$+1)-dimensional case, where $n$ refers to the space and $1$ refers to the time dimension and introduce a geometric formulation of the finite volume method for the general balance laws. We then consider the finite volume approximation on coordinates  to work on the particular  (1+1)-dimensional spacetime for analyzing numerically the relativistic Burgers Equation \eqref{maineq}.

In the following, the introduction of a geometric formulation of the finite volume method on a curved spacetime and its particular case for the local coordinates are given.

\subsection{Geometric Spacetime Finite Volume Method}

Following the article  \cite{ALO}, we  search for a formulation of the finite volume scheme for a hyperbolic balance law given by
\be
\label{GFM}
div  \big(T(v)\big)= S(v),
\ee
posed on an $(n+1)$-dimensional curved spacetime $M$, where  the unknown function $v$ is a scalar field, $ div (\cdot) $ is  the divergence operator, $T(v)$ is the flux vector field and $S(v)$ is the scalar field.

We are now ready to introduce the geometric formulation of the finite method for discretizing the balance laws \eqref{GFM} on $M$.
First  of all, we establish a general triangulation of the spacetimes by
$$\Tcal^h = \bigcup_{K\in \Tcal^h} K
$$
for the manifold $M$ such that $M$ is composed of these spacetime elements $K\subset M$ satisfying the following conditions:
%
\begin{itemize}
\item The boundary  of an element $K$, denoted by $\del K$, is   given by
$$\del K = \bigcup_{e\subset \del K} e,$$ which is piecewise smooth and  includes exactly two ``spacelike'' faces (having an induced Riemannian type metric)
denoted by  {$\ekp$} and {$\ekm$}, and ``timelike''   elements (having an induced Lorentzian type metric) denoted by
$$
e^0 \in \del^0 K:= \del K \setminus \big\{\ekp, \ekm\big\}.
$$
\item The intersection $ K \cap  K'$ of two distinct elements $ K$ and $ K'$ in $\Tcal^h$
is  either a common face of $ K,  K'$ or  is a smooth submanifold with dimension at most $(n-1)$.
\item $|K|$, $|\ekp|$, $|\ekm|$, $|e^0|$ represent the
measures of $K$, $\ekp, \ekm,e^0$, respectively.
\item Along the timelike faces $\ekp,\ekm$, we introduce the outgoing unit normal vector field denoted by $n_{\ekp},n_{\ekm}$, respectively.
\end{itemize}

Furthermore, we assume that $M$ permits a foliation  $\{H_t\}_{t\in(0,\infty)}$,  by oriented spacelike hypersurfaces
such that the parameter $t : M \cup \del M\to [0, +\infty)$ provides us with a global time function. This allows us to separate between future--oriented ($t$ increasing) and past--oriented
($t$ decreasing) timelike directions on $M$.  Thus, $M$ is assumed to be foliated  by hypersurfaces $ {H_t}$ (which are spacelike elements) such that
\be
\label{CAP}
M \cup \del M= \bigcup_{t \geq 0}  {H_t},
\ee
with the time function $t$ if it is determined from a sequence of discrete times $t_0< t_1< t_2< \cdots< t_N$ such that all spacelike faces are submanifolds of ${H_t}=H_{t_n}, (n=0,1, \cdots, N)$ and   $H_0$ is an initial slice.  It follows that  the  class of nonlinear hyperbolic Equations (\ref{GFM}) and (\ref{CAP}) gives a scalar model on which we can analyze numerical methods.

We are ready to  define the finite volume approximations by  averaging the balanced  law (\ref{GFM}) over each element $K \in  \Tcal^h$ of the triangulation. By integrating in space and time, we can write
$$
\int_K {div} (T(v)) =  \int_K S(v).
$$
Applying the Stokes formula to the above equation, it follows that
$$
\int_{\ekp}T^0(v)\, (n_{\ekp},\cdot)  =  \int_{\ekm}  T^0(v)\,  (n_{\ekm},\cdot)
 - \sum_{e^0  \in \del^0 K}\int_{e^0} T^1(v)\,(n_{e^0},\cdot) +
 \int_{ K} S(v).
$$
Using the given averages $v^-_K$ and $v^-_{K_{e^0}} $ computed along $\ekm$ and $e^0\in \del^0 K$,
we need to compute the average $v^+_K$ along $\ekp$. To this aim, we introduce
\be
\aligned
\int_{\ekm}  T^0(v)\, (n_{\ekp},\cdot) &\simeq  |\ekm|  \overline{T}_{\ekm}(v_K^-),\\
\int_{e^0} T^1(v)\,(n_{e^0},\cdot) &\simeq |e^0|\,  Q_{K,e^0}( v_K^-, v_{K_{e^0}}^-),
\endaligned
\ee
and
$$
\int_{ K} S(v)\,  \simeq |K|   \overline{S}_{ K},
$$
where, for each $K$ and $e^0 \in \del ^0 K$, we selected a numerical flux
$Q_{K,e^0} : \RR^2 \to \RR$ satisfying natural properties of consistency, conservation, and monotonicity.
Therefore, the finite volume method of interest takes the form
\be
\label{SHM}
|\ekp|  \,\overline{T}_{\ekp}(v_K^+)=
|\ekm| \, \overline{T}_{\ekm}(v_K^-) - \sum_{e^0\in \del ^0 K} |e^0|  \,Q_{K,e^0} ( v_K^-,v_{K_{e^0}}^-)
+    |K|\,  \, \overline{S}_{ K}.
\ee

Notice that a standard Courant–Friedrichs–Lewy (CFL) condition is assumed for the sake of stability of this scheme. \mbox{For further
details}, the reader may consult  \cite{ALO, DAF, RL, LO}.


\subsection{Finite Volume Schemes in Coordinates}

Following the geometric formulation of the finite volume scheme on $(n+1)$-dimensional spacetime in the previous part, we now take into account the particular case (for $n=1$) and assume that the spacetime is described in coordinates $(t,r)$.
We divide into equally spaced cells $I_j= [ r_{j-1/2}, r_{j+1/2}]$ of size $\Delta r$, centered at $r_j$,
that is,
$$
r_{j+1/2}= r_{j-1/2}  + \Delta r,
$$
satisfying  $$r_{j-1/2} = j \Delta r,\quad r_j= (j+1/2)\Delta r.
$$

Furthermore, we denote by $\Delta t$ the constant time length and we set  $t_n= n \Delta t$.

In order to introduce the finite volume method in local coordinates, we start by rewriting the  hyperbolic balance law \eqref{GFM}  in $(1+1)$ dimension, that is,
\be
\label{new.bal}
\aligned
\del_t T^0(t,r) + \del_r T^1(t,r)=  S(t,r),
\endaligned
\ee
where $T^0, T^1$ are flux fields, and $S$ is the source term.
We integrate Equation \eqref{new.bal} over each grid cell $$[t_n, t_{n+1}] \times [r_{j-1/2}, r_{j+1/2}]$$  to get
$$
\aligned
 \int_{r_{j-1/2}}^ {r_{j+1/2}} &(T^0(t_{n+1},r)- T^0(t_n,r)\,dr
 +  \int_{t_n}^ {t_{n+1}} (T^1(t,r_{j+1/2})- T^1(t,r_{j-1/2}))\,dt\\
 &=\int_{  [t_n, t_{n+1}] \times [r_{j-1/2}, r_{j+1/2}] } S(t,r) \,dt \,dr,
\endaligned
$$
which is equivalent to
$$
\aligned
 \int_{r_{j-1/2}}^ {r_{j+1/2}} T^0(t_{n+1},r)\,dr=& \int_{r_{j-1/2}}^ {r_{j+1/2}} T^0(t_{n},r)\,dr
 -  \int_{t_n}^ {t_{n+1}} (T^1(t,r_{j+1/2})- T^1(t,r_{j-1/2}))\,dt\\
 &+ \int_{  [t_n, t_{n+1}] \times [r_{j-1/2}, r_{j+1/2}] } S(t,r) \,dt \,dr.
\endaligned
$$

The numerical flux functions are approximated by
$$
\aligned
&\widetilde{T}_j ^n\approx {1 \over \Delta r}\int_ {r_{j-1/2}}^ {r_{j+1/2}}T^0(t_{n},r)\, dr,\\
&\widetilde{Q}_{j\pm1/2}^n \approx {1 \over  \Delta t}\int_{t_n}^ {t_{n+1}} T^1(t,r_{j\pm1/2})\,dt,\\
\endaligned
$$
$$
\aligned
&\widetilde{S}_j^n \approx {1 \over {\Delta r \, \Delta t}}  \int_{  [t_n, t_{n+1}] \times [r_{j-1/2}, r_{j+1/2}] } S(t,r) \,dt \,dr.
\endaligned
$$
Using these approximations in the above integrations,  the scheme becomes
\be
\label{version}
\aligned
\widetilde{T}_j ^{n+1}&= \widetilde{T}_j ^n  - \frac{\Delta t}{\Delta r} \big( \widetilde{Q}_{j+1/2}^n -  \widetilde{Q}_{j-1/2}^n   \big)
+ \Delta t\widetilde{S}_j^n.
\endaligned
\ee
Since  $ \widetilde{T}$ is  convex (see \cite{LMO}), by taking the inverse of the above relation, the scheme \eqref{version} is rewritten by
$$
\aligned
v_j ^{n+1}&= \widetilde{T}^{-1}\Big( \widetilde{T}(v_j ^n)  - \frac{\Delta t}{\Delta r} \big( \widetilde{Q}_{j+1/2}^n -  \widetilde{Q}_{j-1/2}^n   \big)
+ \Delta t\widetilde{S}_j^n \Big).
\endaligned
$$

For further details of triangulations and discretization of finite volume approximations  on a  curved spacetime,
we refer the reader to the papers by LeFloch \textit{et al.}  \cite{ALO,  LMO, LM, LO}.

\section{Second Order Godunov--Type Scheme}

In this section, we construct a second order Godunov type finite volume scheme for the derived model which is the main objective of this paper.
Fundamentally, Godunov type algorithms could be summarized by the following steps:
\begin{itemize}
\item Discretize the cells in finite volumes;
\item Obtain the profile reconstruction on the cells;
\item Do specification of the Riemann problems at cell interfaces;
\item Find solution for the Riemann problems;
\item Perform computation for fluxes on cell interfaces;
\item Determine the variables for the next time step.
\end {itemize}

The second order scheme differs from the original Godunov type scheme in reconstruction on the profile and specification of the Riemann problems.
More generally, for the equation
$$
\del_t v + \del_r f(v,r) = 0,
$$
a first order scheme can be converted to a second order method by proceeding the
cell-boundary values which are used in the numerical flux functions to decide the intermediate time level \linebreak $t^{n + 1/2}=(t^n
+t^{n+1})/2$.  In other words, the second order Godunov scheme is
obtained from the edge-boundary values of the reconstructed profile proceeded
by a half time step.

 The construction of the second order type of scheme is formulated by the following steps:
\begin{itemize}
\item Firstly, the variable within the computational cells is reconstructed which gives couples of values $(v_{i,L}^n, v_{i,R}^n)$ in each
computational cell where the letters $L$ and $R$ refer to left and right sides, respectively. Note that  $v_{i,L}^n$ lies between
$v_{i-1}^n$ and $v_{i}^n,$ whereas $v_{i,R}^n$ lies between $v_{i}^n$
and $v_{i+1}^n$.
\item Secondly, one should advance the solution by half a step in time. The intermediate values $v$ at the cell edges
at the time $t^{n + 1/2}=(t^n+t^{n+1})/2$ are denoted by
$(v_{i,L}^n,v_{i,R}^n)$. These values can be calculated by the following:
$$
v_{i,L}^{n+1/2}=v_{i,L}^n-{\Delta t \over 2\Delta
r}[f(v_{i,R}^n)-f(v_{i,L}^n)],
$$
\vspace{-12pt}
$$
v_{i,R}^{n+1/2}=v_{i,R}^n-{\Delta t \over 2\Delta
r}[f(v_{i,R}^n)-f(v_{i,L}^n)].
$$
\item As a third step, we solve  the Riemann problem formed by the intermediate values $(v_{i,L}^n,v_{i,R}^n)$.
The solution $v_{i+1/2}^{n+1/2}$ is used in order to compute the flux
$f_{i+1/2}^{n+1/2}=f(v_{i+1/2}^{n+1/2})$.
\item In the final step, we  advance the solution by the time step $\Delta t$ from $t^n$ using the usual formula
$$
v_i^{n+1}=v_i^n-{\Delta t \over \Delta
r}[f_{i+1/2}^{n+1/2}-f_{i-1/2}^{n+1/2}].$$
\end{itemize}

We refer to the  article by Van Leer  \cite{L} and the book by Guinot \cite{VG} for {construction} of the algorithm
to the second order Godunov type finite volume method.

\section{Numerical Experiments}

\subsection{Implementation of the Godunov Scheme for the Derived Model}

In this part,  numerical tests are displayed  for the model
derived on a   dS  spacetime based on a second order  Godunov
scheme. We prepare the second order finite volume approximation for our
model based on the construction in the previous section.  The method for
the intermediate values and the proceeding solutions are formulated
as follows.
We analyze the given model  with a single shock and
rarefaction for an  initial function considering the Godunov scheme
with a local Riemann problem for each grid cell depending on two
particular cases of cosmological constant $\Lambda$. Note that, for the numerical experiments,
we fix  $r \in [0,1]$. We take into consideration the
fastest wave at each grid cell since
both shock and rarefaction waves are produced in the Riemann problem. We compel transmissive  boundary
conditions on the scheme. By normalization (taking $c=1$) in
Equation \eqref{maineq}, we attain the following model:
\be
\label{godunov1} \del_t v +(1-\Lambda r^2)\del_r\Big({v^2\over
2}\Big)=-\Big(\Lambda r (v^2-1)\Big),
\ee
and the finite volume schemes for this model can be given as
\be
\label{w1}
v_{j\pm
1/2}^{n+1/2}=v_{j\pm 1/2}^n-{\Delta t \over 2 \Delta
r}(b_{j+1/2}^{n} g _{j+1/2}^{n}  - b_{j-1/2}^{n} g _{j-1/2}^{n}
)+{\Delta t \over 2} S_{j\pm 1/2}^{n},
\ee
\be
\label{w2}
v_j^{n+1}=v_j^n-{\Delta t \over \Delta r}(b_{j+1/2}^{n+1/2} g
_{j+1/2}^{n+1/2}  - b_{j-1/2}^{n+1/2} g _{j-1/2}^{n+1/2} )+\Delta t
S_j^{n+1/2}, \ee
with
$$
t^{n + 1/2}=(t^n +t^{n+1})/2
$$
and
$$
\aligned
b_{j\pm1/2}^{n+1/2}&=1-\Lambda({r}_{j\pm1/2}^{n+1/2})^2,\\
g _{j-1/2}^{n+1/2}&=f(v_{j-1}^{n+1/2},v_{j}^{n+1/2}),\\
g_{j+1/2}^{n+1/2}&=f(v_j^{n+1/2},v_{j+1}^{n+1/2}),\\
\endaligned
$$
where the source term is
$$
S_j^{n+1/2}=-\Big(\Lambda r_j^{n+1/2}((v_j^{n+1/2})^2-1)\Big).
$$
Moreover, the flux function $f(v_1,v_2)$ is  written as follows
\be \label{flux}
f(v_1,v_2)= \begin{cases}{v_1^2 \over 2}, \quad  \, \mbox {if} \qquad v_1>v_2 \qquad   \mbox {and} \qquad v_1+v_2>0\\
{v_2^2 \over 2}, \quad  \, \mbox {if} \qquad v_1>v_2 \qquad   \mbox {and} \qquad v_1+v_2<0\\
{v_1^2 \over 2}, \quad  \, \mbox {if} \qquad v_1\leq v_2 \qquad   \mbox {and} \qquad v_1>0\\
{v_2^2 \over 2}, \quad  \, \mbox {if} \qquad v_1\leq v_2 \qquad   \mbox {and} \qquad v_2<0\\
\,{0},\quad  \, \, \, \mbox{if}  \qquad v_1\leq v_2 \qquad \mbox{and} \qquad  v_1 \leq 0 \leq v_2\\
              \end{cases}.
\ee
Notice that  the speed term is
\be
\label{spd} 1-\Lambda({r}_{j\pm1/2}^{n+1/2})^2,
\ee and, in order for the stability condition in the scheme to be satisfied,
we choose $\Delta t$ and $\Delta r$ so that the CFL condition satisfied, \textit{i.e.},
\be
\label{CFL}
{\Delta t \over \Delta r}
\underset{j}{\max}\Big|1-\Lambda({r}_{j\pm1/2}^{n+1/2})^2\Big| \leq
1.
\ee

The implementation of the second order Godunov scheme is based on
the construction given above. Observe that $\Lambda>0$  and $\Lambda=0$ refer to the dS and Minkowski geometries, respectively. We compare the particular cases
$\Lambda=1>0$ and $\Lambda=0$ for the schemes of the
model with shock and rarefaction waves; that is, we analyze the model via the second order Godunov scheme for two different geometries.

\begin{figure}
\label{a22} \centering
\includegraphics[width=5cm]{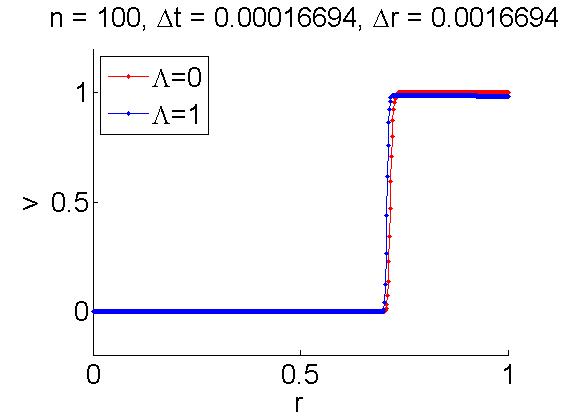} \includegraphics[width=5cm]{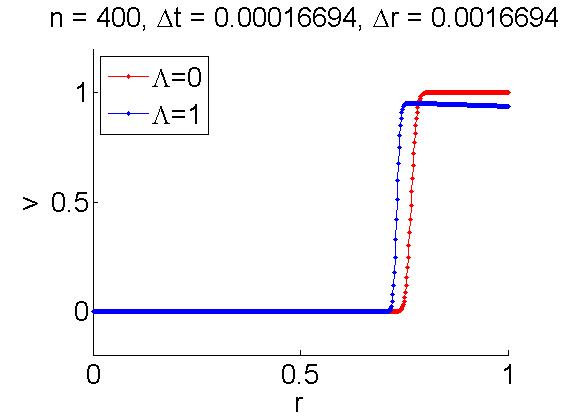}
\includegraphics[width=5cm]{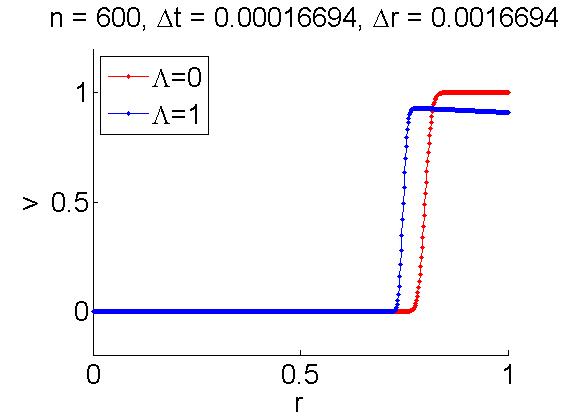} \includegraphics[width=5cm]{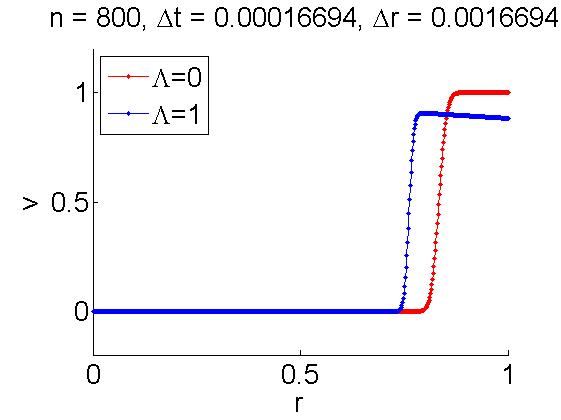}
\caption{The numerical solutions given by the second order Godunov scheme with a
rarefaction for $\Lambda=0$ and $\Lambda=1$.}
\end{figure}

\begin{figure}
\label{a23} \centering
\includegraphics[width=5cm]{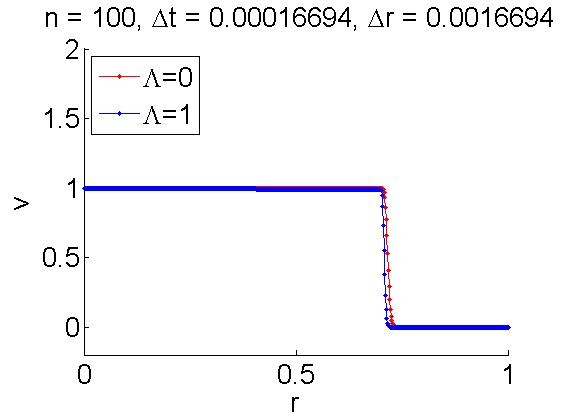} \includegraphics[width=5cm]{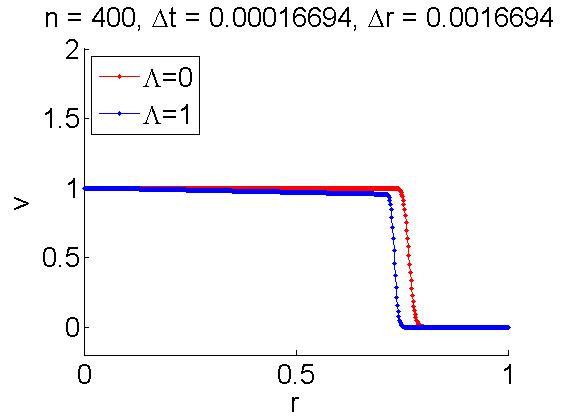}
\includegraphics[width=5cm]{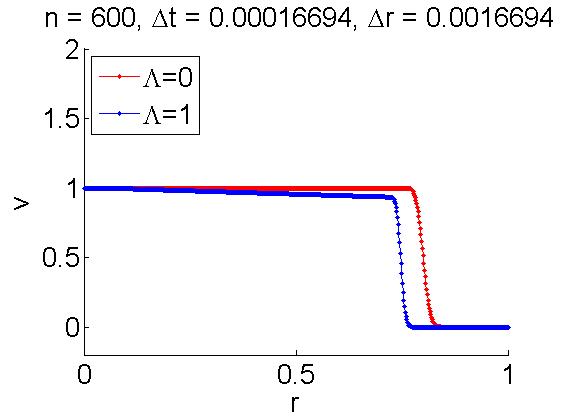} \includegraphics[width=5cm]{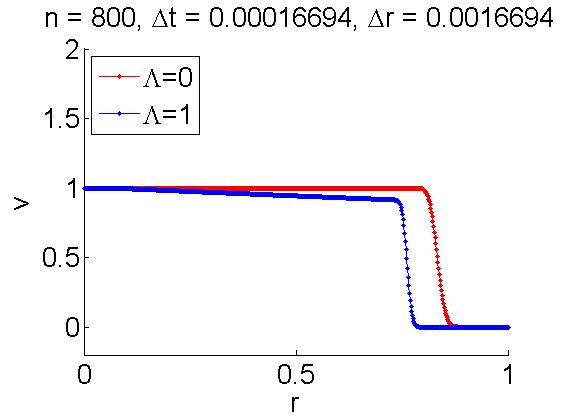}
\caption{The numerical solutions given by the second order Godunov scheme with a
shock for $\Lambda=0$ and $\Lambda=1$.}
\end{figure}

 The results are sketched in Figures \ref{a22} and \ref{a23}. The rarefactions are drawn in Figure \ref{a22} and shocks are given in Figure \ref{a23} for two  particular cases of $\Lambda$. The iterations are performed for different numbers of iterations; for instance, for the first figure of Figure $1$-$2$, the number of iterations $n=100$, for the second figure $n=400$, for the third one $n=600$, and for the last one $n=800$.  Moreover, the CFL condition \eqref{CFL} is strictly less than $1$ in each test for stability of the scheme.
We observe that the numerical solution for the case $\Lambda=0$ (represented by the red curve)  moves faster than the case $\Lambda=1$ (represented by the blue curve). In other words, we have a relatively higher  speed when the spacetime is flat comparing with the dS geometry.
This fact can also be verified numerically by substituting  $\Lambda=0 $ and $  \Lambda=1$ into the speed term
given by \eqref{spd}, \textit{i.e.}, by the term
$$
1-\Lambda({r}_{j\pm1/2}^{n+1/2})^2,
$$
 which feeds the theoretical background of
the model with the numerical results.
%
%


%

We also observe that the negative values of $\Lambda$ yield the velocity to become greater than $1$  by Equation \eqref{spd}, which is not acceptable since they exceed the speed of light (as we assumed that $c=1$ for the numerical calculations). In other words, the sign of the cosmological constant $\Lambda$ is significant for the numerical results by which the geometry under consideration is identified.


%

\section{Conclusions}

A new nonlinear hyperbolic model,
which describes the propagation and interactions of shock waves on
the dS and flat spacetimes, is derived in this  study. We take into account  the relativistic
Euler system on a curved background and impose a zero
pressure in the statement of the energy--momentum tensor for
perfect fluids in order to derive our model \eqref{maineq}.
 This model involves a cosmological
constant   $\Lambda$, which can be normalized by taking the
 positive, negative  and zero values.  According to the positive, zero and negative  values of $\Lambda$, the geometry of interest becomes the de Sitter, Minkowski and Anti-de Sitter spacetimes, respectively. We have the following remarkable observations:

\begin{itemize}
\item  The standard Burgers equation is recovered
when $\Lambda$ is chosen to vanish;
\item  Various mathematical properties concerning the hyperbolicity,
genuine nonlinearity, shock waves, and rarefaction waves are established;
\item Investigation for the class of static and spatially homogenous solutions is carried out;
\item Shock wave solutions to the model for the different possible values of the cosmological constant $\Lambda$ are investigated.
\end {itemize}

The analysis for two
possible values of the cosmological constant $\Lambda$ relies on a numerical scheme,
which applies discontinuous solutions and is based on the finite volume technique. We establish a second-order Godunov scheme for our model which gives the main originality to the current paper.

\begin{itemize}
\item The scheme is consistent with the conservative form of  our model, which yields correct computations of  weak solutions containing shock waves.
\item The numerical solution curves for the Anti-de Sitter spacetimes yields the velocity to become greater than one, which is not acceptable as it exceeds the speed of light.
\item One of the most obvious findings emerging from this study is that it allows us to make a numerical comparison of the relativistic models on the de Sitter and Anti-de Sitter spacetimes by means of finite volume approximations. Indeed, the method is not efficient for the relativistic model on the Anti-de Sitter spacetimes, whereas it is efficient for the de Sitter spacetime.
\item The efficiency and convergency of the numerical method depends on the value of the cosmological constant $\Lambda$. The sign of  $\Lambda$ is significant. If $\Lambda$ is smaller than zero, the method is not efficient and the solution curves do not converge. On the other hand, if $\Lambda=0$ and $\Lambda>0$, the method gives efficient and convergent numerical results.
\item The numerical experiments illustrate the convergence, efficiency and robustness of the proposed scheme on the de Sitter background.
\item Numerical solutions for the cases  $\Lambda=0$  and $\Lambda=1>0$ are compared. The solution curve corresponding to $\Lambda=0$ moves {\bf faster} than the solution curve corresponding to $\Lambda=1>0$   (Figures $1$ and $2$). This feature can be clarified by observing that the characteristic speed
$$
 1-\Lambda({r}_{j\pm1/2}^{n+1/2})^2
 $$
is upgraded
by reduced $\Lambda$  in the Equation \eqref{maineq}.
\end {itemize}

To sum up, we highlight that the existence of the model of the concerning equation that permits one to improve and analyze numerical methods for shock capturing schemes and to attain definite conclusions concerning the convergence, robustness and efficiency of the schemes. According to the  background geometry, different techniques may be needed to guarantee that certain classes of time-dependent or space-dependent solutions for the model {be}
preserved by the scheme.  A strategic perspective can be given as to investigate relativistic Burgers equations on different backgrounds with more complicated points.

%
%
\subsection*{Acknowledgments} The second author
is supported by {\it the Rectorate of Middle East Technical
University (METU) } through the grant {\bf ``Project BAP-01-01-2016-007''}.


\bibliographystyle{mdpi}

\begin{thebibliography}{999} 




\bibitem{CLO}
 {Ceylan, T.; LeFloch, P.G.;  Okutmustur, B}. The relativistic Burgers equation on an FLRW background and its finite volume approximation. \textbf{2015}, arXiv:1512.08142v1.

\bibitem{CO}
 {Ceylan, T.;  Okutmustur, B}. Relativistic Burgers equation on a de Sitter
background. Derivation of the model and finite volume approximations. {\em Int. J. Pure Math.} {\bf 2015}, {\em  2}, 21--29.

\bibitem{CO2}
{Ceylan, T.;  Okutmustur, B}. Finite volume approximation of the relativistic Burgers equation on a Schwarzschild-(anti-)de Sitter spacetime. {\bf 2016}, submitted.

\bibitem{LMO}
{LeFloch, P.G.; Makhlof, H.;  Okutmustur, B.} Relativistic Burgers equations on a curved spacetime. Derivation and finite volume approximation. {\em SIAM J. Num. Anal.} {\bf 2012}, {\em 50},  2136--2158.

\bibitem{Ceylan}
Ceylan, T.  The Relativistic Burgers Equation on an FLRW Background and Its Finite Volume Approximation. { PhD Thesis}, Middle East Technical University,  Ankara,  Turkey, 2015.

\bibitem{LM}
 {LeFloch, P.G.;  Makhlof, H}. A geometry-preserving finite volume method for compressible fluids for Schwarzshild spacetime. {\em Commun. Comput. Phys.} {\bf 2014}, {\em 15}, 827--852.

\bibitem{ALO}
{Amorim, P.;  LeFloch, P.G.; Okutmustur, B}. Finite volume schemes on Lorentzian manifolds. {\em Commun. Math. Sci.} {\bf 2008}, {\em  6},  1059--1086.

\bibitem{SNK}
{Kruzkov, S.N}. First order quasilinear equations in several independent variables. {\em Math. USSR Sb.} {\bf1970}, {\em 81},  228--255.

\bibitem{GST}
 {Groah, J.M.; Smoller, J.; Temple, B.} {\em Shock Wave Interactions in General Relativity}; Springer: New York, NY, USA, 2007.

\bibitem{Gosse}
{Gosse, L.} Locally inertial approximations of balance laws arising in (1+1)-dimensional general relativity. {\em SIAM J. Appl. Math.}  {\bf 2015}, {\em 75},  1301--1328.

\bibitem{VG}
{Guinot, V.} {\em Godunov--type Schemes: An introduction for Engineers}; Elsevier: Amsterdam, The Netherlands, 2003.

\bibitem{L}
{Van Leer, B.} On the relation between the upwind-differencing schemes of Godunov,
Engquist-Osher and Roe. {\em SIAM J. Sci. Stat. Comput.} {\bf 1984},  {\em 5}, 1--20.

\bibitem{W}
 {Wald, R.M}. {\em General Relativity}; The University of Chicago Press: Chicago, IL, USA, 1984.

 \bibitem{DAF}
 {Dafermos, C.M.} {\em Hyperbolic Conservation Laws in Continuum Physics}; Springer: Berlin, Germany, 2000.

 \bibitem{RL}
  {LeVeque, R.J.}  {\em Finite Volume Methods for Hyperbolic Problems}; Cambridge Texts in Applied Mathematics; Cambridge University Press: Cambridge, UK, 2002.
%


\bibitem{LO}  {LeFloch P.G.;  Okutmustur B}.
{Hyperbolic conservation} laws on spacetimes. A finite volume scheme based on differential forms. {\em Far. East J. Math. Sci.} {\bf 2008}, {\em  31},  49--83. 
\end{thebibliography}
\makeatletter
\renewcommand\@biblabel[1]{#1. }
\makeatother

\end{document}